# Computing Nilpotent Quotients in Finitely Presented Lie Rings

Csaba Schneider
Kossuth Lajos University of Arts and Sciences
Department of Algebra and Number Theory
4010 Debrecen pf.12 (Hungary)
`schcs@jordan.math.klte.hu`

April 20, 1996

**Abstract**

A nilpotent quotient algorithm for finitely presented Lie rings over $\mathbb{Z}$ (LieNQ) is described. The paper studies the graded and non-graded cases separately. The algorithm computes the so-called nilpotent presentation for a finitely presented, nilpotent Lie ring. A nilpotent presentation consists of generators for the abelian group and the products expressed as linear combinations for pairs formed by generators. Using that presentation the word problem is decidable in $L$. Provided that the Lie ring $L$ is graded, it is possible to determine the canonical presentation for a lower central factor of $L$. LieNQ's Complexity is studied and it is shown that optimizing the presentation is NP-hard. Computational details are provided with examples, timing and some structure theorems obtained from computations. Implementation in C and GAP interface are available.

# 1  Introduction

The nilpotent quotient algorithm for finitely presented Lie rings—LieNQ as we shall refer to it in the sequel—operates with Lie rings over $\mathbb{Z}$. Its essential goal is to compute a sufficient and useful presentation for a given Lie ring. It actually means to compute a presentation for its abelian group and to compute the structure constants in order to determine the Lie ring structure. By efficient and useful we mean that several important pieces of information can be read off immediately (e.g.: nilpotency, nilpotency class etc). Furthermore, by this presentation the so-called *word problem* is decidable, although it is known to be undecidable in the most general case. on word problem in groups and Lie algebras see (Hall, 1958).

Such algorithms have existed for several decades for groups and Lie rings as well. Recall, that groups and Lie rings have similar structure if we think the commutator as a second operation in a group. It has a Jacobi-like property and an identity similar to the distributivity and it is anti-commutative, see for example (Robinson, 1982). This allows us to alter the known group algorithms for Lie rings.

The first widely known (and maybe the most successful) algorithm is due to (Havas & Newman, 1980). It has several implementation in FORTRAN, C and GAP. The description of these were reported in (Celler, Newman, Nickel & Niemeyer, 1993) and (Havas & Newman, 1980). They were successfully used to compute the order of some Burnside groups.

From our point of view, probably the most important algorithms are presented in (Havas, Newman & Vaughan-Lee, 1990), (Nickel, 1993) and (Nickel, 1995a). The algorithm described in this paper is a mixture of the latter mentioned ones, however it differs at several points from those.

There are three main sections in this paper. The first part contains the basic properties of the nilpotent presentation and describes an algorithm to compute it in the most general case. The next part is devoted to graded Lie rings. We shall see that in the graded case we can simplify LieNQ and, using other techniques, it is possible to determine the isomorphism type of the underlying abelian group of such a Lie ring. The last part contains computational issues. We shall see how the implementation works and some examples with timing will display the power of LieNQ in both graded and non-graded cases.

The basic properties of Lie rings can be found in any textbook on that topic. We refer (Humphreys, 1972) for details. The first lemma is of fundamental importance for our goal, since it essentially claims that it is realistic to build a nilpotent quotient algorithm for finitely presented Lie rings.

**Lemma 1.1** *In a finitely presented Lie ring each lower central factor is finitely presented.*

The proof can be done by a similar argument to that presented in (Robinson, 1982) 5.2.6, but one might want to see (Schneider, 1996) for more details.

## 2 The Nilpotent Presentation of a Lie Ring

Let $L$ denote a Lie ring and $L^i$ the $i$th term of its lower central series. Product in $L$ is denoted by brackets—as it is usual—and we use the left-normed convention that is, $[a, b, c] = [[a, b], c]$. If $L$ is finitely presented Lemma 1.1 implies that $L^i/L^{i+1}$ is finitely presented for each $i$. Since at the present stage we are only interested in the nilpotent factor $L/L^k$, we think $L$ is nilpotent of class $k - 1$. In other words we think $L/L^k$ is equal to $L$. In this case $L$ has a so-called *nilpotent presentation*, that is a tuple $(\mathcal{G}, S)$, where $\mathcal{G}$ is a finite set of generators, $\{a_1, \ldots, a_r\}$ say, and $S$ is a finite set of relators of the form

$$\gamma_i \cdot a_i = \alpha_{i,i+1} \cdot a_{i+1} + \cdots + \alpha_{i,r} \cdot a_r \quad \text{for some } i, \qquad \text{TR}$$

$$[a_j, a_i] = \beta_{j,i,j+1} \cdot a_{j+1} + \cdots + \beta_{j,i,r} \cdot a_r \quad \text{for each } j > i, \qquad \text{PR}$$

where the $\alpha$, $\beta$, $\gamma$ are integer coefficients, furthermore, $\gamma_i > 1$ where applicable. Note that the relations of type TR are sometimes referred to as *torsion relations*, while those of PR are known as *product relations*. A more detailed description and the proof for existence of such presentation can be found in (Schneider, 1996), while it is easy to see, that a nilpotent presentation always defines a nilpotent Lie ring. Observe that we demand the existence of TR relations for each pair $(j, i)$, where $j > i$, but only some of those of type PR (even possibly none). Throughout this paper $I$ denotes the set of indices, so that $i \in I$ iff $a_i$ has TR relation.

One easily sees that the relations of type TR determine the underlying abelian group structure for $L$, while the relations of type PR define its ring structure. The $\beta$ can be viewed as structure constants for the $\mathbb{Z}$-module defined by $\mathcal{G}$ and the TR relations.

For our purpose this concept will be too general, so we shall keep some restrictions. Define a weight function $\omega : \mathcal{G} \to \mathbb{N}$ as follows. Let $\omega$ be an increasing function in $i$, such that the following holds.

1. $\omega(a_1) = 1$,

2. if the nilpotent presentation contains a relation of the form

   $$[a_j, a_i] = w_{ji},$$

   where $w_{ji}$ is a sum of integer multiples of generators from $\mathcal{G}$, then all generators $a_k \in \mathcal{G}$ such that $\omega(a_k) < \omega(a_j) + \omega(a_i)$ have coefficients zero in $w_{ji}$.

An other restriction, if $a_k \in \mathcal{G}$ and $\omega(a_k) > 1$, we require $a_k$ to have a definition, i.e. a relation of the form

(1) $$[a_j, a_i] = a_k$$

where $\omega(a_k) = \omega(a_j) + \omega(a_i)$ and $\omega(a_i) = 1$. It might easily happen that a nilpotent presentation contains more then one relations of the form (1) for some $i$, in this case we arbitrarily choose one for the definition for $a_i$.

From now on by a nilpotent presentation we shall always mean weighted nilpotent presentation which satisfies the condition stated above for the existence of definitions.

Suppose that we are given a nilpotent presentation

(2) $$\langle\, \mathcal{G} \,|\, \gamma_i \cdot a_i = w_i,\ i \in I, [a_j, a_i] = w_{ji},\ 1 \leq i < j \leq r \,\rangle.$$

Let $L$ denote the Lie ring presented by (2). Clearly (2) does not care about Jacobi identities, so $L$ is not automatically coincides with the freest non-associative anty-commutative ring holding the relations of (2) in it. However we can arrange a collection process in $L$. Suppose, we are given a word-sum $\ell$ expressed in term of the elements from $\mathcal{G}$. The collection process consists of two main steps. In the first step we substitute the products by the right hand side of the appropriate PR relation. In the second step we reduce those coefficients that are not less than the corresponding $\gamma_i$ in (2).

After proceeding so, we get a normal form for $\ell$ subject to (2). That is,

(3) $$\ell = \alpha_1 \cdot a_1 + \cdots + \alpha_r \cdot a_r$$

where $0 \leq \alpha_i < \gamma_i$ whenever $i \in I$. An arbitrary $\ell$ might have more than one normal form of type (3). However uniqueness of (3) is a crucial point, if we want to attack the word problem in $L$. We might want the nilpotent presentation to possess this property according to the following definition.

**Definition 2.1** *A nilpotent presentation is said to be consistent if every element of $L$ uniquely has a normal form.*

In Section 3.3 we shall see how to check if an arbitrary nilpotent presentation is consistent and develop a method to make it consistent if it is not. In general, it is an equivalent condition to that of Definition 2.1 that the element 0 uniquely has a normal form.

From the computational point of view it is important to see that all the relations of type PR are not necessary to determine the structure of $L$. Our first result is devoted to this observation.

**Lemma 2.1** *Let $L$ be a Lie ring given by the nilpotent presentation*

$$\langle\, \mathcal{G} \,|\, \gamma_i \cdot a_i = w_i,\ i \in I, [a_j, a_i] = w_{ji},\ 1 \leq i < j \leq r \,\rangle.$$

*Then $L$ is already determined by the presentation*

(4) $$\langle\, \mathcal{G} \,|\, \gamma_i \cdot a_i = w_i,\ i \in I,\ [a_j, a_i] = w_{ji},\ 1 \leq i < j \leq r,\ \omega(a_i) = 1 \,\rangle.$$

*Proof.* Essentially we want to prove that all PR relations can be expressed by those listed in (4). We proceed by an induction argument on $\omega(a_i)$. If $\omega(a_i) = 1$ then we are done. Suppose now, that $\omega(a_i) = k > 1$. Using the definition for $a_i$ and the identities of a Lie ring one has: $[a_j, a_i] = [a_j, [a_k, a_l]] = -[a_k, a_l, a_j] = [a_l, a_j, a_k] + [a_j, a_k, a_l] = [a_j, a_k, a_l] - [a_j, a_l, a_k]$. And the latter ones are known by the hypothesis. In the first equality we substituted $a_i$ by its definition. $\square$

# 3 Computing the Nilpotent Presentation

The aim of this section is to show how it is possible to compute a consistent nilpotent presentation for a nilpotent factor of a a given finitely presented Lie ring $L$ and determine an epimorphism from $L$ onto that nilpotent factor.

The basic ideas—besides being simple—are similar to those used in (Havas, Newman & Vaughan-Lee, 1990) and (Nickel, 1995a). In fact, it is based on an induction argument. Suppose, we are given a consistent nilpotent presentation for $L/L^l$ and an epimorphism from $L$ onto $L/L^l$, we shall show how to extend them to a consistent presentation and an epimorphism for $L/L^{l+1}$.

As computational tools, matrices over $\mathbb{Z}$ play an important rôle. Computing the *row Hermite normal form* is a crucial point of the algorithm. Since its details can be found in (Nickel, 1995a) (Schneider, 1996), (Sims, 1994), we do not emphasize them further in this paper.

## 3.1 Computing the Abelian Factor

Suppose that we are given a Lie ring $L$ with finite presentation $\langle\, \mathcal{X} \mid \mathcal{R} \,\rangle$, where $\mathcal{X}$ is a finite set of generators, while $\mathcal{R}$ is a finite set of relators. For simplicity, suppose $\mathcal{X} = \{x_1, \ldots, x_n\}$. Set $\mathcal{G} = \{a_1, \ldots, a_n\}$. Construct a map $\phi : \mathcal{X} \to \mathcal{G}$, by simply saying
$$x_i \phi = g_i.$$
Furthermore we set all PR relations to be trivial and $I = \emptyset$. We call the presentation constructed above a *trivial nilpotent presentation* and it is clear that it represents the free abelian Lie ring $\mathcal{F}_n^{ab}$ on $n$ generators. It is straightforward to see that it contains $L/L^2$ as a factor ring. Indeed, let $\mathcal{I}$ be the ideal generated by the set $\mathcal{R}\phi$, where we evaluate $\phi$ as if it was a homomorphism elementwise on $\mathcal{R}$. Then we have

**Lemma 3.1** $L/L^2$ *is isomorphic to* $\mathcal{F}_n^{ab}/\mathcal{I}$.

*Proof.* Both Lie rings are generated by $n$ elements and satisfy exactly the same relations. It checks the isomorphic property. $\square$

What we do in practice is the following. Evaluate the relations in $\mathcal{F}_n^{ab}$ as described above and put the images in a matrix $M$. Then compute the row Hermite normal form $M^H$ for $M$. It is well known that the rows of $M$ generate the same subgroup, in the free abelian group, as those of $M^H$. If the first non-zero element, the $i$th say, of a row of $M^H$ is 1, i.e. we have a row of the form
$$(0, \ldots, 0, 1, m_{i+1}, \ldots, m_n)$$
we simply throw $a_i$ out of $\mathcal{G}$ and modify the map $\phi$
$$x_i \phi = -m_{i+1} \cdot a_{i+1} - \cdots - m_n \cdot a_n.$$
Otherwise, if the leading element is greater then one, i.e. we have a row of the form
$$(0, \ldots, 0, m_i > 1, m_{i+1}, \ldots, m_n),$$

we introduce a TR relation

$$m_i \cdot a_i = -m_{i+1} \cdot a_{i+1} - \cdots - m_n \cdot a_n.$$

and put $i$ into $I$.

Proceeding each row as described above, we obtain a presentation for $L/L^2$. We set the surviving generators of weight one and it is straightforward to see that the presentation obtained so far is indeed a presentation for the abelian factor and $\phi$ extends to an epimorphism. We summarize this result in

**Proposition 3.1** *The so obtained presentation is a consistent nilpotent presentation for $L/L^2$ and $\phi$ extends to a epimorphism from $L$ onto $L/L^2$.*

The untouched images under epimorphism will be referred to as definitions of generators of weight one. This terminology will be important later on.

## 3.2 Extending the Presentation

Now we shall see how the induction argument works in details. Suppose, we have a nilpotent presentation for the factor ring $L/L^l$, an epimorphism from $L$ onto $L/L^l$. Extending the presentation will consist of three steps, i.e.

1. extending the epimorphism $\phi$,
2. modifying the TR relations,
3. modifying the PR relations.

Extending the epimorphism means introducing a new generator for all $x_i \in \mathcal{X}$, where $x_i\phi$ is not a definition. In other words, if $x_i\phi = w_i$, we modify $x_i\phi$ so that $x_i\phi = w_i + t_i$.

The torsion relations will change so that if $i \in I$, $\gamma_i \cdot a_i = w_{ii}$, then we alter $\gamma_i \cdot a_i = w_{ii} + t_{ii}$.

Those of the PR relations that are not definitions are modified in a similar way. If $[a_j, a_i] = w_{ji}$ $j > i$ and $[a_j, a_i]$ is not a definition, let the new relation be $[a_j, a_i] = w_{ji} + t_{ji}$.

Note, that all newly introduced generators $t_i, t_{ii}, t_{ji}$ are different from one another. We introduce PR relations so that all of them are central. At this stage we do not alter $I$. We prove the following

**Proposition 3.2** *The so extended presentation contains $L/L^{l+1}$ as a factor ring.*

*Proof.* Since we did not alter the definitions, every relation holding in the extended presentation is satisfied in $L/L^{l+1}$. □

In practice we use Lemma 2.1. This lemma essentially says, that we do not need to introduce new generators for all PR relations, but only for those of the form $[a_j, a_i] = w_{ji}$ for $j > i$ and $a_i$ is of weight one. The remaining ones can be computed as we have seen in the proof of this lemma. This computation is usually referred to as *computing the tails*.

## 3.3 Enforcing Consistency

Recall, we called a nilpotent presentation consistent if each element in $L$ uniquely had a normal form. The extended presentation investigated in the previous subsection might happen to be not consistent, as the case is in general.

In order to investigate the uniqueness property, we introduce an operation on the nilpotent quotient. First of all let $A$ be an abelian group generated by $\mathcal{G}$ and the TR relations of the extended presentation. Introduce a map $\psi$ such that,

$$\psi : \mathcal{G} \times \mathcal{G} \to A \qquad \psi(a_j, a_i) = \begin{cases} w_{ji} & \text{if } j > i, \\ -w_{ji} & \text{if } j < i, \\ 0 & \text{if j=i} \end{cases},$$

provided that the nilpotent presentation possesses a PR relation of the form $[a_j, a_i] = w_{ji}$ for $j > i$. Then the following is true.

**Proposition 3.3** $\psi$ can be extended to a binary operation on $A$ iff

$$o_j \cdot [a_j, a_i] = \sum_{k=j+1}^{n} \alpha_{jk}[a_k, a_i], \qquad C1$$

where $j \in I$ and the relation corresponding to $a_j$ is of the form

$$o_j \cdot a_j = \alpha_{j,j+1} a_{j+1} + \cdots + \alpha_{in} \cdot a_n.$$

*Proof.* It is immediate that the condition is necessary. To see that it is sufficient as well, observe that the two expressions, $w$ and $w'$ say, are equal, it means that they can be transformed to each other by using TR relations. The condition of the proposition claims that such a transition respects the equality $\psi(w, v) = \psi(w', v)$ for any $v$. Similar result is obtained in the second variable of $\psi$ by swapping the two variables. □

If $A$ is equipped by that operation, it happens to be a non-associative ring. We prove the following

**Lemma 3.2** *The nilpotent presentation of $L$ is consistent iff $A$ is a Lie ring.*

*Proof.* Observe, $L = A/J$ where $J$ is the ideal generated by all instances of Jacobi identity. A normal word $\ell$ represents the zero element in $L$ iff $\ell \in J$. If $A$ is a Lie ring then $J = 0$, so $\ell$ must be the trivial word. On the other hand, if the nilpotent presentation is consistent, then all instances of Jacobi identity collect to the trivial word, so $J = 0$ and $A = L$ holds. □

One can easily check the Jacobi identity, since it is a multi-linear identity, thus it is enough to check it on the triples formed by the abelian group generators. However we can restrict ourselves even more, by the following

**Lemma 3.3** *A nilpotent presentation $(\mathcal{G}, S)$ of a Lie ring $L$ is consistent iff the Jacobi identity holds for each triple $(a_i, a_j, a_k) \in \mathcal{G} \times \mathcal{G} \times \mathcal{G}$ where $1 \leq i < j < k \leq n$ and $a_i$ is of weight one. In other words*

$$[a_i, a_j, a_k] + [a_j, a_k, a_i] + [a_k, a_i, a_j] = 0 \qquad C2$$

*for $1 \leq i < j < k \leq n$ and $a_i$ is of weight one.*

The above lemma appeared in its first form in (Vaughan-Lee, 1984) for pc-presentation of $p$-groups. It was modified for Lie algebras in (Havas, Newman & Vaughan-Lee, 1990), where one finds its proof. The Lie ring case can be handled in the same way as Lie algebras.

### 3.4  Enforcing the Defining Relations

One can easily observe that the presentation obtained so far is a presentation of the freest Lie ring of class $l$ that happens to contain $L/L^{l+1}$ as a factor ring in it. What is still needed is to assure that the images of the elements of $\mathcal{R}$, under the epimorphism, are 0. It simply means taking the factor ring over the ideal generated by the epimorphic images. By the induction hypothesis those images vanish in $L/L^l$, so they must lie in $L^i/L^{i+l}$ and we are allowed to use abelian group methods again, such as Hermite normal form and other integer based strategies.

In practice we proceed as follows. Put all elements obtained from the consistency relations (C1), (C2) and the defining relations together in a matrix $M$. Compute the row Hermite normal form $M^H$ for $M$. Using those rows of $M^H$ whose leading element is equal to 1, we can eliminate some of the generators from the nilpotent presentation. The other rows express linear dependence among the generators over $\mathbb{Z}$. They can be viewed as TR relations and added to the presentation. Extend the weight function to the newly introduced generators, by saying that they are of weight $i$ and add the the indices of generators that were provided with TR relations to $I$. The following lemma makes sure that generators of weight $i$ have definitions.

**Lemma 3.4** *If $L$ is a Lie ring and suppose that $L/L^2$ is additively generated as*

$$L/L^2 = \langle \ell_1 + L^2, \ldots, \ell_k + L^2 \rangle$$

*and $L^{i-1}/L^i$ is additively generated as*

$$L^{i-1}/L^i = \langle \ell_1' + L^i, \ldots, \ell_l' + L^i \rangle$$

*then $L^i/L^{i+1}$ is additively generated as*

$$L^i/L^{i+1} = \langle [\ell_m', \ell_n] + L^{i+1} \mid 1 \leq m \leq l,\, 1 \leq n \leq k \rangle.$$

The proof of the above lemma can be done essentially in the same way as in the group case. That proof can be found in (Sims, 1994) Proposition 2.6. But see also (Schneider, 1996) for a detailed study of the Lie ring case.

# 4 The Graded Case

## 4.1 Simplifying LieNQ in Graded Lie Rings

The following section is devoted to graded Lie rings. Gradedness of a Lie ring is defined in accordance with the following definition.

**Definition 4.1** *Let $L$ be a Lie ring, $L$ is said to be graded if $L$ splits into a direct sum of its lower central factors.*

Since a Lie ring of that kind possesses a relatively easy structure, we are allowed to simplify our algorithm in this case.

First of all, we need not extend the TR relations any more since if $a_i \in L^l/L^{l+1}$ for some $i \in I$ implies $\gamma_i \cdot a_i \in L^l/L^{l+1}$, so the right hand side of the relation corresponding to $a_i$ cannot contain generators of weight greater than $l$.

Similar thing happens when we extend the PR relations. At the $l$th step we only alter the relations of the form

$$[a_j, a_i] = w_{ji} \qquad \text{for } i > j \quad i + j = l.$$

We regard these facts when we introduce new generators and compute the tails.

Now it should be clear that enforcing the consistency simply means enforcing the consistency relations (C2) for the triples $(a_i, a_j, a_k) \in \mathcal{G} \times \mathcal{G} \times \mathcal{G}$, where $1 \leq i < j < k \leq n$ $\omega(a_i) + \omega(a_j) + \omega(a_k) = l$ and $\omega(a_i) = 1$. And the relations of type (C1) are checked for only pairs $(a_i, a_j) \in \mathcal{G} \times \mathcal{G}$, where $\omega(a_i) + \omega(a_j) = l$. Note that $l$ still denotes the nilpotency class of the current factor.

We are allowed to simplify it further, namely we need to enforce only those relations of $\mathcal{R}$ which have weight $l$.

The natural question arises. How can one recognize a graded Lie ring by looking at its finite presentation? We only mention here the well known fact that a Lie ring defined by homogeneous relations is always graded.

## 4.2 Another approach—Canonical Nilpotent Presentation

In a graded Lie ring we shall construct the so-called *canonical nilpotent presentation* instead of nilpotent presentation defined in Section 2. The concept of the canonical nilpotent presentation fits the concept of canonical presentation for finitely generated abelian groups. In the sequel we shall use additive notation for abelian groups and + to denote the binary operation in such structures. Recall that $C_n$ denotes the cyclic group of order $n$.

**Theorem 4.1** *Let $(A, +)$ finitely generated abelian group. Then $A$ can be written as a direct sum of cyclic groups, i.e.*

(5) $$A \cong C_{k_1} \oplus \cdots \oplus C_{k_r} \oplus C_\infty \oplus \cdots \oplus C_\infty,$$

*where $k_1|k_2,\ldots,k_{r-1}|k_r$. The (5) form is called the canonical decomposition for $A$. The canonical decomposition is unique for an arbitrary such $A$.*

The proof of the theorem can be found in (Suzuki, 1982) Theorem 5.2.

In the following we define the canonical nilpotent presentation for a graded Lie ring.

**Definition 4.2** *A Lie ring presentation is said to be canonical nilpotent presentation, if it is of the form*

(6) $\qquad \langle a_1, \ldots a_n \,|\, c_i \cdot a_i \quad \text{for } i \in I,\ [a_j, a_i] = w_{ij} \quad \text{for } j > i \rangle,$

*where the $w_{ij}$ are normal words as in (3).*

**Remark 4.1** *Using the results of Section 3 we shall see that each finitely presented graded nilpotent Lie ring possesses such a presentation. Essentially it is enough to see, that the generators corresponding to a certain lower central factor can be chosen so that the torsion relations are trivial in (6). Conversely, it is rather easy to see that a Lie ring defined by a canonical nilpotent presentation is always a finitely presented, nilpotent and graded.*

**Remark 4.2** *Actually we shall construct more than a simple presentation defined in (6). We shall construct a presentation in which the lower central factors are canonically presented as abelian groups, according to Theorem 4.1. That is, we shall enforce the appropriate divisibility conditions for the $c_i$.*

The following theorem is an easy consequence of Theorem 4.1 and Definition 4.2

**Theorem 4.2** *Requiring the divisibility condition for each lower central factor in (6), the coefficients $\gamma_i$ are uniquely determined by the isomorphism type of $L$.*

### 4.3 Integer Matrices

In this section we collect the background knowledge needed to construct our nilpotent quotient algorithm. First of all, a matrix is called *integer matrix* if its entries are rational integers.

**Definition 4.3** *Let $A$ and $B$ $m \times n$ be integer matrices. $A$ and $B$ are called equivalent if there exist $P$ $m \times m$ and $Q$ $n \times n$ integer unimodular matrices so that $PAQ = B$.*

**Remark 4.3** *The unimodular matrix $P$ corresponds to elementary row operations, while $Q$ corresponds to elementary column operations. An elementary row (column) operation is negating a row (column), adding an integer multiple of a row (column) to an other row (column), swapping two rows (columns).*

**Remark 4.4** *A trivial but rather useful observation is that the Smith normal form of the matrix $(a_1, \ldots a_n)$ is $(\gcd(a_1, \ldots, a_n), 0 \ldots, 0)$.*

The following result can be traced back to (Smith, 1861).

**Theorem 4.3** *For each M $m \times n$ integer matrix there exists an S $m \times n$ diagonal matrix, such that the non-zero entries of S are non-negatives and $s_{i-1,i-1}|s_{i,i}$ for $2 \leq i \leq \min(m,n)$. S is called the Smith normal form for M.*

**Remark 4.5** *Theorem 4.1 and Theorem 4.3 look very similar. Indeed, they are closely related. It is not hard to see if we put the defining relators of an abelian group into a matrix, then abelian groups with equivalent relation matrices are isomorphic. Theorem 4.1 has a standard proof as follows. Put the defining relations of A into a matrix M. Then compute the unique Smith normal form S for M and read off the canonical presentation for A from S.*

(Hartley & Hawkes, 1970) provides a naive algorithm to compute Smith normal form for a given integer matrix. However, in practice it does not seem to be widely applicable. Besides the fact, that the intermediate entries are growing very fast, we shall need to compute the unimodular matrix $Q$ too. $Q$ is not unique and we might want to choose a capable $Q$ with moderate entries. The above mentioned algorithm does not fit this goal as examples will point it out later.

Werner Nickel (Nickel, 1995b) suggested the following algorithm.

**Algorithm 1**
INPUT: $m \times n$ *matrix M*
OUTPUT: *Smith normal form for M*
**while** *M is not diagonal*
    **begin**
        $M := $ *Hermite normal form of M*
        $M := $ *Transpose of M*
    **endwhile**
*Enforce divisibility on M's diagonal*
**return** $M$
**end**

**Lemma 4.1** *Algorithm 1 computes the Smith normal form in finitely many steps.*

*Proof.* If after a transposition the first entry of the $i$th column divides all entries in the its column, we can reduce that column without producing trash in the upper half.

If it does not divide some entry in its column, after a reduction the first element will become a smaller positive integer. After finitely many steps it will divide all column entries (in the worst case it becomes 1). □

More sophisticated algorithms exist to compute Smith normal form as reported in (Havas & Majewski, 1994), (Havas & Majewski, 1995) and (Havas, Holt & Rees, 1993). They use various techniques such as LLL, modular methods and pivoting strategy. They behave very well in practice, however, our matrices are rather sparse and Algorithm 1 is suitable for our purposes.

**Example 4.1** *In (Sims, 1994) on Page 379 one finds an example with detailed analysis of its Smith normal form computation. Using the naive algorithm it was hopeless to compute Smith normal form for the matrix. Using Algorithm 1 one can compute Smith normal form, producing 5276851227 as the maximum entry of Q.*

*Algorithms with more sophisticated background do not behave much better. Its Hermite normal form was computed using ideas described in (Havas, Holt & Rees, 1993), implemented by B.Majewski in* GAP. *Norm driven pivoting strategy produced* $34522312156188926219838487089136 89 \sim 3.4 \cdot 10^{34}$ *as maximum entry in the transformer matrix. Their LLL based algorithm computed the Smith normal form with 298901 as maximum entry of P.*

**Example 4.2** *The other example arose from the following Lie ring presentation.*

$$\langle x, y \,|\, [y, x, x, y], [y, x, x, x] + [y, x, y, y] \rangle$$

*Computing its 10th factor,* LIENQ *produced an* $30 \times 16$ *matrix. Naive approach implemented in* C, *with 32 bit arithmetic, produced integer overflow. Algorithm 1 provided a transformer matrix with largest entry 12. Norm driven Hermite normal form computation produced 62, LLL based algorithm produced 5 as largest entry in P.*

The idea behind Algorithm 1 is that we keep the number of column operations as low as possible. It seems to be a good idea using LLL based algorithm, however, it has worse complexity than that of traditional strategies.

As we saw the main problem is to keep magnitude of the entries of the intermediate matrices and the transforming matrix under control. To find the optimal strategy—even in the simplest case—is hopeless, as stated in (Havas & Majewski, 1994).

**Proposition 4.1** *Let us given* $a = (a_1, \ldots, a_n)$, *a vector consisting of n positive integers. To find a vector x, such that x is shortest either in* $L_0$ *or in* $L_\infty$ *norm and solves the equality* $x \cdot a = \gcd(a1, \ldots, a_n)$ *is* NP-*hard.*

The following proposition is an easy modification of Theorem 3 in (Havas & Majewski, 1994).

**Proposition 4.2** *Let us be given an integer vector* $a = (a_1, \ldots, a_n)$. *The task of minimizing the vectors in the unimodular matrix Q, where* $aQ = (gcd(a_1, \ldots, a_n), 0 \ldots, 0)$, *with respect to the maximum norm (*MINGCD*), is* NP-*hard.*

For the definition and basic properties of NP-hard problems we refer (Davis & Weyuker, 1983).

Let the abelian group $A = \langle x_1, \ldots x_n \,|\, r_1, \ldots, r_m \rangle$ be given by generators and relations. As it was mentioned above, the canonical presentation can be constructed for $A$ by computing the Smith normal form for the relation matrix of $A$. We emphasize this idea in details now.

The relation matrix is constructed as follows. If $r_i = a_{i,1}x_1 + \cdots + a_{i,n}x_n$ let $M = (a_{i,j})$ the $m \times n$ matrix. Suppose that the Smith normal form for M is

$$S = \begin{pmatrix} s_{11} & \cdots & 0 & \cdots & 0 \\ \vdots & \ddots & \vdots & & \\ 0 & \cdots & s_{rr} & & \vdots \\ \vdots & & & \ddots & \\ 0 & & \cdots & & 0 \end{pmatrix}.$$

Then

(7) $$A \cong C_{s_{11}} \oplus C_{s_{22}} \oplus \cdots \oplus C_{s_{rr}} \oplus C_\infty \oplus \cdots \oplus C_\infty,$$

where the number of infinite factors is $n - r$. Since $S = PMQ$ for suitable unimodular matrices $P$ and $Q$, $S$ can be obtained from $M$ by elementary row and column operations. Since a row operation does not alter the generating set of $A$, $P$ can be forgotten, while $Q$ can be used to obtain the isomorphism between the two presentations of $A$.

**Proposition 4.3** *Let $A$, $M$, $Q$ and $S$ be as above. Suppose moreover that the cyclic factors in (7) are generated by $a_1$, $a_2$,...,$a_n$ respectively. Then*

$$\begin{pmatrix} a_1 \\ \vdots \\ a_n \end{pmatrix} = Q^{-1} \begin{pmatrix} x_1 \\ \vdots \\ x_n \end{pmatrix} \quad \text{and} \quad Q \begin{pmatrix} a_1 \\ \vdots \\ a_n \end{pmatrix} = \begin{pmatrix} x_1 \\ \vdots \\ x_n \end{pmatrix}.$$

*Proof.* See (Sims, 1994) Proposition 8.3.1. □

### 4.4 Computing the Canonical Nilpotent Presentation

After this preparation we are ready to describe, how one can compute canonical nilpotent presentation for a finitely presented nilpotent graded Lie ring.

Recall, we computed the matrix consisting of consistency relations and defining relations. We computed its Hermite normal form. Now we shall compute its Smith normal form together with the transformer matrix $Q$ and its inverse $Q^{-1}$. Using Proposition 4.3 one can compute new generators for the canonical presentation. Put, $Q = (q_{i,j})$ and $Q^{-1} = (\hat{q}_{i,j})$.

The condition for the existence of the definition cannot be kept any more in that case. We use the following approach instead.

Suppose, at a certain layer we introduced new generators $x_1, \ldots, x_s$, so that they correspond to the following products:

$$x_i = [b_{1i}, b_{2i}] \quad \text{where} \quad \omega(b_{1i}) = 1 \quad \text{and} \quad \omega(b_{2i}) = \omega(x_1) - 1,\ 1 \leq i \leq s.$$

Using Proposition 4.3 one gets the generators $a_1,\ldots,a_s$ of the canonical presentation for that layer. Moreover one has that

$$\begin{pmatrix} a_1 \\ \vdots \\ a_n \end{pmatrix} = Q^{-1} \begin{pmatrix} [b_{11}, b_{21}] \\ \vdots \\ [b_{1s}, b_{2s}] \end{pmatrix}$$

So, we can introduce definitions for the canonical generators of the following form

$$a_i := \hat{q}_{i,1}[b_{11}, b_{21}] + \cdots + \hat{q}_{i,s}[b_{1s}, b_{2s}] \quad \text{for } 1 \le i \le s.$$

That is, the definitions are linear combinations consisting of Lie products with generators of smaller weight.

In Section 2 we already computed products of current weight. They need one more look as follows. After tail computation, we expressed those products in terms of newly introduced generators $x_1,\ldots,x_s$. Now, those generators are to be replaced by the canonical generators $a_1,\ldots,a_s$. To do that we use Proposition 4.3 again. That is,

$$x_i = q_{i,1}a_i + \cdots + q_{i,s}a_s \quad \text{for } 1 \le i \le s.$$

We proceed the old generators occuring in the epimorphic images in the same way.

Lemma 2.1 is still true in the canonical case and its proof can be done in essentially the same way as presented before. Lemma 3.3 seemingly has no connection to the form of the definitions. However, having a deeper look at its proof, one finds that it heavily depends on the concept of definitions, but can be modified to the new approach too.

The canonical nilpotent presentation computed so far depends on the transformer unimodular matrix $Q$. Among many possible matrices, we want to choose the one with the smallest possible entries, because they will appear in the presentation. This is a hard problem as stated in

**Theorem 4.4** *Let L be a finitely presented nilpotent, graded Lie ring. The task of minimizing the vectors in the right hand side of the epimorphic images and the PR relations, according to the maximum norm, (MINLIE) is NP-hard.*

*Proof.* We reduce MINGCD to MINLIE in polynomial time. Let us be given an instance $(c_1, \ldots, c_n)$ of MINGCD. Write down the finite presentation

$$\langle x_1, \ldots, x_n \mid c_1 x_1 + \cdots + c_n x_n = 0, [x_i, x_j] = 0 \quad \text{for} \quad 1 \le i < j \le n \rangle.$$

Clearly, this can be done in time proportional to $n^2$. Use LieNQ on this presentation. LieNQ terminates at the abelian factor and produces the following output:

(8)
$$\begin{aligned} x_i\varphi &= q_{i,1}a_1 + \cdots q_{i,n}a_n \quad \text{for } 1 \le i \le n, \\ \gcd(c_1, \ldots, c_n)a_1 &= 0, \\ [a_i, a_j] &= 0 \quad \text{for } 1 \le j < i \le n. \end{aligned}$$

Where the transformer matrix $Q = (q_{i,j})$. The ability of minimizing the vectors in (8) implies the ability if minimizing the vectors in $Q$ and hence solving MINGCD. □

An other disadvantage of this approach is presented by the following theorem.

**Theorem 4.5** *The absolute values of the coefficients appearing in the canonical nilpotent presentation cannot be bounded by a function depending only on the number of generators.*

*Proof.* Easy consequence of (Havas & Majewski, 1994) Lemma 4. □

## 5 Performance

### 5.1 Implementing LieNQ in C

The implementation for the LieNQ algorithm has been written in C programming language. The currently available version is the Version 2.0, which provides time measure, computations for graded and non-graded cases separately.

The program takes a finite presentation $(\mathcal{X}, \mathcal{R})$ for a Lie ring $L$ and the nilpotency class of required factor as its input. It outputs the TR, PR relations and the images of the elements of $\mathcal{X}$ under the epimorphism. The nilpotency class can be omitted, in this case the program runs until the lower central series stabilizes or the algorithm goes beyond the capacity of the computer.

Throughout the computation the data of the nilpotent presentation is stored in normal word form. At previous stages they were stored as coefficient vectors, but it turned out to be more efficient in an array consisting of structures, where the first element of a structure contains the number of a generator and the second one its coefficient (the idea is due to (Nickel, 1995b)). For instance, the Lie ring element
$$0 \cdot a_1 + 0 \cdot a_2 + 2 \cdot a_3 - 4 \cdot a_4 + 0 \cdot a_5$$
is stored as
$$\big((3, 2), (4, -4), (0, 0)\big).$$
The first component of the last pair indicates an `EOW` (end of word) sign, while the second one, in fact, is arbitrary.

The information of the nilpotent presentation is stored in several arrays as follows:

| | |
|---|---|
| `Coefficients[i]` | $\gamma_i$ if $i \in I$, 0 otherwise, |
| `Power[i]` | $w_i$, right hand side of the TR relation $\gamma_i \cdot a_i$, |
| `Product[j][i]` | $w_{ji}$, right hand side of the PR relation $[a_j, a_i]$, |
| `Epimorphism[i]` | $x_i\phi$, the $i$th epimorphic image, |
| `Weight[i]` | $\omega(a_i)$, the weight of $a_i$, |
| `Dimension[i]` | number of generators of weight $i$, |
| `Definition[i]` | the definition of $a_i$. |

Note that `Definition[i]` is an array of structures. A component of the array either contains the number of two generators and a positive integer, if the definition of $a_i$ is like $a_i = \sum c_i \cdot [a_{k_i}, a_{l_i}]$, or contains the number of the generator from $\mathcal{X}$ in the first component and 0 in the second and the third if $a_i$ is defined as an epimorphic image. In the latter case the length of the array is, of course, one.

The defining relations are stored as expression trees. See (Nickel, 1995a) for its description which emphasizes its advantages in details.

During the computation graded and non-graded cases are distinguished according to whether the user switched on the `-g` option or not. There are different tails routines, that is `Tails` and `GradedTails`, and consistency routines, like `Consistency` and `GradedConsistency`, built in the code. Of course, introducing new generators is also a different task.

It is, however, important to introduce the new generators in a right order. Lemma {commlemma makes possible to express all tails in terms of generators $a_k$ of weight $l$ where $a_k = [a_j, a_i]$ for some $a_j$ of weight $l-1$ and $a_i$ of weight 1. In order to use this fact in practice, we introduce those generators the last and when we compute the row Hermite normal form for the matrix they will be, in fact, the only surviving generators and all the others disappear. In this way, we ensure that the new generators indeed have definitions, so the restrictions we kept are reasonable. In the graded case such restrictions are not kept.

The main routines for matrix computation and the basics of the parser program are due to W. Nickel. The first mentioned one uses the GNU MP package to deal with long integers.

LieNQ was successfully compiled on several Unix platforms, such that Linux, Free-BSD, SunOS, Solaris.

## 5.2  GAP Interface

LieNQ is might be invoked from GAP (Schönert *et al.*, 1994) using a GAP interface. Since GAP does not know about Lie algebras until version 3.5, it is necessary to have at least GAP-3.5 to use the interface. We remark, that even GAP-3.5 knows only *Lie algebra*s but not Lie rings over a ring, so LieNQ kills torsions in that case and GAP considers the result as Lie algebra over $\mathbb{Q}$.

The interface consists of three functions, two of them are auxiliary ones and one is worth to more interest. An auxiliary function is implemented to print a Lie ring presentation according to LieNQ's syntax, and one is implemented to print a usage-message in the case of possible errors.

Now we describe how to use the function `LieNilpotentQuotient`.

`LieNilpotentQuotient(P, c)`
`LieNilpotentQuotient(P)`

In the first case the function computes the canonical nilpotent presentation for the $c$th nilpotent factor of $P$. In the second case it computes the largest nilpotent factor of $P$ if it exists.

In both cases $P$ has to a be record consisting of at least two fields: generators and relators. A relator is arbitrary Lie element, using the left-normed convention. In complex cases it is a good idea to give relators as strings, since addition is not defined between two lists.

The function returns a table Lie algebra record containing a field for the finite presentation $P$, one for the epimorphism and one for the lower central series of the Lie algebra defined by the output presentation.

## 6   Some Sample Computations

In this section we present some examples, that was tried by LieNQ. We consider the following Lie ring presentations.

$L_1 = \langle x, y \rangle$

$L_2 = \langle x, y \mid [y, x, y], [y, x, x, x, x, x] \rangle$

$L_3 = \langle x, y \mid |[[y, x], x, y], [[y, x], y, x], [[y, x], x, x] + [[y, x], y, y] \rangle$

$L_4 = \langle a, b, c, d, e \mid [b, a], [c, a], [e, c], [e, d], [d, a] = [c, b], [d, b] = [e, a], [d, c] = [e, b] \rangle$

$L_5 = \langle e_1, e_2, e_3, e_4, e_5, e_6, e_7, e_8, e_9, e_{10} \mid$

$[e_2, e_1, e_1], [e_1, e_2, e_2], [e_3, e_1], [e_4, e_1], [e_5, e_1], [e_6, e_1], [e_7, e_1], [e_8, e_1],$

$[e_9, e_1], [e_{10}, e_1], [e_3, e_2, e_2], [e_2, e_3, e_3], [e_4, e_2], [e_5, e_2], [e_6, e_2], [e_7, e_2],$

$[e_8, e_2], [e_9, e_2], [e_{10}, e_2], [e_4, e_3, e_3], [e_3, e_4, e_4], [e_5, e_3], [e_6, e_3], [e_7, e_3],$

$[e_8, e_3], [e_9, e_3], [e_{10}, e_3], [e_5, e_4, e_4], [e_4, e_5, e_5], [e_6, e_4], [e_7, e_4], [e_8, e_4],$

$[e_9, e_4], [e_{10}, e_4], [e_6, e_5, e_5], [e_5, e_6, e_6], [e_7, e_5], [e_8, e_5], [e_9, e_5], [e_{10}, e_5],$

$[e_7, e_6, e_6], [e_6, e_7, e_7], [e_8, e_6], [e_9, e_6], [e_{10}, e_6], [e_8, e_7, e_7], [e_7, e_8, e_8],$

$[e_9, e_7], [e_{10}, e_7], [e_9, e_8, e_8], [e_8, e_9, e_9], [e_{10}, e_8], [e_{10}, e_9, e_9], [e_9, e_{10}, e_{10}],$

$[e_3, e_2, e_1, e_2], [e_4, e_3, e_2, e_3], [e_5, e_4, e_3, e_4], [e_6, e_5, e_4, e_5], [e_7, e_6, e_5, e_6],$

$[e_8, e_7, e_6, e_7], [e_9, e_8, e_7, e_8], [e_{10}, e_9, e_8, e_9] \rangle$

The first example presents the free Lie ring of rank two. The second and third examples are found in (Caranti, Mattarei, Newman & Scoppola, 1994) and are of great importance from the point of view of thin Lie algebras. The fourth example is due to Jürgen Wisliceny and has a very nice structure, as we shall see later. The fifth one is not else than the presentation of the positive part of the classical Lie algebra $A_{10}$ factored out by the ideal generated by the torsions appearing over $\mathbb{Z}$.

The structure of $L_4$ is very nice as we mentioned before. Having a look at its canonical nilpotent presentation one has

**Theorem 6.1** *The underlying abelian group of the lower central factors of $L_4$ is of the form*

$$L^{i+1}/L^{i+2} \cong \begin{cases} C_{i/2} \oplus C_\infty \oplus C_\infty \oplus C_\infty & \text{if } i \text{ is even,} \\ C_\infty \oplus C_\infty \oplus C_\infty \oplus C_\infty \oplus C_\infty & \text{if } i \text{ is odd.} \end{cases}$$

*for* $0 \leq i \leq 67$.

By computation, the following structure theorem is true, for $L_5$.

**Theorem 6.2** *$L_5$ is nilpotent of class 10. Its lower central factors are of the form*

(9) $$L_5^i/L_5^{i+1} = C_\infty \oplus \cdots \oplus C_\infty,$$

*where the number if direct factors is $11 - i$ for $1 \leq i \leq 10$.*

The structure of $L_1$ is well known. (Hall, Jr., 1959) Theorem 11.2.2 provides a formula to compute the dimension for a lower central factor of $L_1$. The result of the computation coincides with the known dimensions.

The following table contains the computational informations concerning those examples.

| example | class computed | nr. of gens. | torsion rank | CPU time (ms) |
|---------|----------------|--------------|--------------|---------------|
| $L_1$   | 10             | 226          | 0            | 811100        |
| $L_2$   | 14             | 89           | 70           | 66200         |
| $L_3$   | 12             | 81           | 63           | 33150         |
| $L_4$   | 7              | 51           | 22           | 41333         |
| $L_5$   | 6              | 45           | 0            | 254016        |

Since all Lie rings presented by those relations are graded we used the graded approach to compute the presentation for them. That is summarized in the following table.

| example | class computed | nr. of gens. | torsion rank | CPU time (ms) |
|---------|----------------|--------------|--------------|---------------|
| $L_1$   | 10             | 226          | 0            | 8633          |
| $L_1$   | 12             | 757          | 0            | 257383        |
| $L_2$   | 14             | 57           | 38           | 3150          |
| $L_2$   | 18             | 195          | 171          | 197483        |
| $L_3$   | 12             | 81           | 63           | 4133          |
| $L_4$   | 67             | 300          | 32           | 903553        |
| $L_5$   | 10             | 55           | 0            | 36083         |

All computations, except for the 67th class of $L_4$, were done on a 486-DX2 100 MHz PC with 16 MB memory plus 16 MB swap space. The above mentioned exception was computed on a Sun Sparccenter 2000 with 256 MB of RAM.

One easily sees that using those ideas speeds up the computation with a significant factor. The need of using better algorithm for Smith normal form computation is reasonable. While computing $L_3$ an integer overflow occured with 32 bit arithmetic at the 68th factor.

# References

Andrea E. Caranti and Sandro Mattarei and Mike F. Newman and Carlo M. Scoppola (1994), "Thin Groups of Prime-power Order and Thin Lie Algebras". Manuscript.

Frank Celler and M.F. Newman and Werner Nickel and Alice C. Niemeyer (1993), "An algorithm for computing quotients of prime-power order for finitely-presented groups and its implementation in GAP", Research Report **SMS-127-93**, School of Mathematical Sciences, Australian National University.

Martin D. Davis and Elaine J. Weyuker (1983), *Computability, Complexity and Languages*. Academic Press.

Philip Hall (1958), "Some Word Problems", *Journal of London Math. Soc.*, **33**, 482–496.

Marshall Hall, Jr. (1959), *The Theory of Groups*. Macmillan Co., New York.

B. Hartley and T.O. Hawkes (1970), *Rings, Modules and Linear Algebra*. Chapman and Hall, London.

George Havas and Derek F. Holt and Sarah Rees (1993), "Recognizing Badly Presented $Z$-Modules", *Linear Algebra Appl.*, **192**, 137–163.

George Havas and Bohdan S. Majewski (1994), "Hermite Normal form Computation for Integer Matrices", *Congressus Numerantium*, **105**, 87–96.

George Havas and Bohdan S. Majewski (1995), "Integer Matrix Diagonalization", *J. Symbolic Comput.*, **11**.

George Havas and M.F. Newman (1980), "Application of computers to questions like those of Burnside", *Burnside Groups*, Lecture Notes in Math., **806**, (Bielefeld, 1977), pp. 211–230. Springer-Verlag, Berlin, Heidelberg, New York.

George Havas and M.F. Newman and M.R. Vaughan-Lee (1990), "A nilpotent quotient algorithm for graded Lie Rings", *J. Symbolic Comput.*, **9**, 653–664.

Humphreys, J.E (1972), *Introduction to Lie Algebras and Representation Theory*. Graduate Text in Mathematics. Springer-Verlag New York, Heidelberg, Berlin.

Werner Nickel (1993), *Central extensions of polycyclic groups*, PhD thesis. Australian National University.

Werner Nickel (1995a), "Computing Nilpotent Quotients in Finitely Presented Groups", *Geometric and Computational Prospectives of Finite Groups*, DIMACS series.

Werner Nickel (1995b), "Personal Communication".

Derek J. Robinson (1982), *A Course in the Theory of Groups*, Graduate Texts in Math., **80**. Springer-Verlag, New York, Heidelberg, Berlin.

Csaba Schneider (1996), *Computations in Finitely Presented Lie Rings*, Master's thesis. Kossuth Lajos University of Arts and Sciences.

Martin Schönert *et al.* (1994), GAP – *Groups, Algorithms and Programming*. Lehrstuhl D für Mathematik, RWTH, Aachen.

Charles C. Sims (1994), *Computation with finitely presented groups*. Cambridge University Press.

Henry J.S. Smith (1861), "On Systems of Linear Indeterminate Equations and Congruences", *Philos. Trans. Royal Soc. London*, **cli**, 293–326.

Michio Suzuki (1982), *Group Theory I*, Grundlehren Math. Wiss., **247**. Springer-Verlag, Berlin, Heidelberg, New York.

M.R. Vaughan-Lee (1984), "An Aspect of the Nilpotent Quotient Algorithm", *Computational Group Theory*, (Durham, 1982), pp. 76–83. Academic Press, London, New York.